\documentclass{article}
\input{preamble} % All my macros

\title{Near-the-Axis Universality in Last Passage Percolation}
\author{Sam McKeown \and Xinyi Zhang} 

\begin{document}
\maketitle

\begin{abstract}
    This note establishes a universal directed landscape limit for last passage percolation models in an intermediate scaling regime. We find as a quick consequence the transversal fluctuations for geodesics taken near the axis. We extend the technique of Bodineau and Martin, who in \cite{bodi-mart} have already shown universal one-point fluctuations in this regime. 
\end{abstract}

\section{Introduction}
\label{sec:intro}
The directed landscape \cite{dov} is the conjectured universal scaling limit for a wide class of random planar metrics and other models belonging to the Kardar-Parisi-Zhang (KPZ) class in $1 + 1$ dimensions. Last passage percolation (LPP) is one such model, wherein we consider i.i.d weights $\set{\omega(a) : a \in \bZ^2}$ attached to the vertices of the integer lattice $\bZ^2$. Between vertices $a,\, b \in \bZ^2$, with $a \le b$ componentwise, these weights define a \emph{passage time}
\begin{equation}
\label{eq:LPPDef}
    L_\omega(a; b) = \max_{\pi \in \Pi(a; b)}\sum_{r \in \pi} \omega(r).
\end{equation}
Here $\Pi(a; b)$ consists of all weakly increasing nearest-neighbour paths on the lattice $\bZ^2$ between $a$ and $b$, namely those of the form $\pi = (\pi_0, \pi_1, \dots, \pi_n)$ with $\pi_{k + 1} - \pi_k \in \set{e_1, e_2}$ and $\pi_0 = a$, $\pi_n \in \set{b - e_1, b - e_2}$. We use $e_1$ and $e_2$ to denote the horizontal and vertical unit vectors, respectively.

The conjectured joint convergence of the passage times $\set{L_\omega(a; b) : a,\, b \in \bZ^2}$ to the directed landscape has only been rigorously established for geometric and exponential weights. A closely related model for which this limit is also known is Brownian last passage percolation (BLPP). Instead of the lattice $\bZ^2$, our model now lives on the semi-continuous space $\bR \times \bZ$, and the weights are replaced by a collection of independent two-sided Brownian motions $\set{B_k(x) : k \in \bZ,\, x \in \bR}$. Passage times in this setting are defined analogously: if $(x, m) \le (y, n)$, then we define
\begin{equation}
\label{eq:BLPPTimes}
    L_B(x, m; y, n) = \max_{x = z_{m - 1} \le z_{m} \le \cdots \le z_{n} = y}\sum_{i = m}^{n} \sqrbrac{B_{i}(z_{i}) - B_{i}(z_{i - 1})}.
\end{equation}

The desired directed landscape limit looks at rescalings of the variables $\set{L_\omega(a; b)}$ when $a,\, b$ are inside the box $[0, n] \times [0, n]$, then takes $n \to \infty$. Such limits remain out of reach for general weight distributions. However, Bodineau and Martin considered in \cite{bodi-mart} a different scaling using thin rectangles $[0, n] \times [0, n^\alpha]$, with $\alpha < 3 / 7$, and showed that the Koml\`os-Major-Tusn\`ady (KMT) coupling of the weights to Brownian motion can be used to transfer knowledge of BLPP to general LPP in this regime. Specifically, they show that under the appropriate normalisation, the passage times $L_\omega\brac[\big]{0, 0; n, n^\alpha}$ converge to the Tracy-Widom GUE distribution. A result of this type was obtained independently by Baik and Suidan in \cite{baik-suidan}, using instead Skorokhod embedding to relate the environment to Brownian motion.

Neither pair of authors extracted richer information from this comparison, owing to the fact that much of the technology regarding the scaling limits of planar directed metrics did not yet exist. In this note we observe that modern techniques, namely the tools developed by Dauvergne and Virag in \cite{dauv-vir-lis} and improvements on the KMT coupling due to Sakhanenko, can extend this idea beyond one-point convergence to a full directed landscape limit. The universal transversal fluctuations of geodesics conjectured in \cite{bodi-mart} fall out as a consequence.

\subsubsection*{Acknowledgements}
% I had previously pu "we" here -- sorry to put words in your mouth!
The first author thanks Dominik Schmid for suggesting this as a tractable and more interesting alternative to the problem which the first author had initially mentioned.

\section{Preliminaries}
\label{sec:prelim}
We briefly introduce the central objects in the study of such scaling limits. We do not attempt to precisely define these or consider any of the attendant technicalities, and refer the reader to \cite{dov,dauv-vir-lis} for a more thorough treatment. 

From here on, we normalise the weights so that $\Ex{\omega} = 0$ and $\Exabs{\omega}^2 = 1$. Also assume that we have some $p > 2$ such that $\Exabs{\omega}^p < \infty$. Notice that we do not assume the distribution to be bounded above or below. This means that, after reversing signs, everything we say for LPP will also hold in directed first passage percolation.

\subsection{Melons and the Airy line ensemble}
By a \emph{line ensemble} we will mean a (possibly finite) sequence of functions $F = (F_1, F_2, \dots)$, with each $F_k : \bR_{\ge0} \to \bR$ continuous. We will also see line ensembles with domain $\bR$.

An important transformation on such line ensembles is the \emph{melon map}, a version of the Robinson–Schensted–Knuth correspondence with many pleasing properties \cite{dauv-nica-vir-rsk}. Given an integer $n \ge 1$ and points $(x, m),\, (y, n) \in \bR_{\ge0} \times \bN$, with $x \le y$ and $1 \le m \le n - k + 1$, we define $Q^k(x, m; y, n)$ to be the set of $k$-tuples of cadlag step functions $\pi = (\pi^1, \dots, \pi^k)$ such that:
\begin{itemize}
    \item Each $\pi^i : [x, y] \to \set{m, \dots, n}$ is a non-decreasing cadlag function.
    \item We have $\pi^i(y) = n$ for all $i$.
    \item If $i > j$, then $\pi^i (u) > \pi^j (u)$ for all $u \in [x, y)$.
\end{itemize}
We may think of such a $\pi$ as a collection of \emph{directed}, \emph{non-intersecting paths} on $\bR \times \bN$ with their endpoints at $(y, n)$.

Define integrals with respect to such paths by 
\begin{equation*}
    \int d F \circ \pi^i = \sum_{j = m}^{n} F_j(x^i_{j + 1}) - F_j({x^i_{j}}),
\end{equation*}
where the times $x^i_j$ are defined by $[x^i_j, x^i_{j + 1}) = \set{z : \pi^i(z) = j}$, or $x^i_j = x^i_{j + 1}$ if this set is empty. All of this is just to say that the integral of an ensemble along a path $\pi^i$ is the increment of the function $F_{\pi^i(z)}(z)$ on the interval $[x, y]$, less the jumps. An integral with respect to a collection $\pi = (\pi^1, \dots, \pi^k)$ is just 
\begin{equation*}
    \int d F \circ \pi = \sum_{i = 1}^{k} \int d F \circ \pi^i.
\end{equation*}
 
We now define the melon $W^n F = ((W^n F)_1, \dots, (W^n F)_{n})$. The components are determined by
\begin{equation}
\label{eq:MelonDef}
    \sum_{i = 1}^{k}(W^n F)_i(y) = \sup_{\pi \in Q^k(0, 1; y, n)} \int d F \circ \pi.
\end{equation}
In particular, $(W^n F)_1(y)$ has the form of a last passage time as defined in \eqref{eq:BLPPTimes}, from $(0, 1)$ to $(y, n)$. Observe that melons are non-intersecting: for all $y \ge 0$,
\begin{equation*}
    (W^n F)_1(y) \ge (W^n F)_2(y) \ge \cdots \ge (W^n F)_n(y).
\end{equation*}

The \emph{Airy line ensemble} is an object conjectured to be a universal scaling limit for such melons. It is a random non-intersecting line ensemble $\frakA = (\frakA_1, \frakA_2, \dots)$, where the $\frakA_k : \bR \to \bR$ are almost surely continuous and satisfy
\begin{equation*}
    \frakA_1(y) > \frakA_2(y) > \cdots
\end{equation*}
for all $y \in \bR$. The next theorem will provide our working definition of the Airy line ensemble.

\begin{theorem}
\label{thm:BLPPALELimit}
    Let $B = (B_1, B_2, \dots)$ consist of standard Brownian motions. The Airy line ensemble can be found as the limit
    \begin{equation*}
        n^{1 / 6}\brac[\big]{(W^n B)_k (1 + 2 y n^{-1 / 3}) - 2 \sqrt{n} - 2 y n^{1 / 6}} \distto \frakA_k(y),
    \end{equation*}
    the convergence here being in distribution, uniform on compact sets.
\end{theorem}

We take the above statement from Theorem 2.4 of \cite{dov}, but this form of convergence was first proved by Johansson in \cite{johan-png-airy}. It follows from looking at the top lines of our ensembles that $\frakA_1(y)$ serves as the scaling limit for $L_B(0, 1; 1 + 2 y n^{-1 / 3}, n)$.

A more intrinsic characterisation of the Airy line ensemble was proven, with much effort, by the authors of \cite{agg-huang-ale}. It is the unique random line ensemble such that $\frakA_1(y) - y^2$ is stationary, and which satisfies the so-called Brownian Gibbs property. Essentially, the distribution of $\frakA$ restricted to a compact set, conditioned on the boundary values, is equal to that of non-intersecting Brownian bridges.

\subsection{The Airy sheet and directed landscape}
The Airy line ensemble provides the scaling limit for a one-parameter slice of BLPP, but of course we would like to take a limit in all four parameters jointly. The \emph{Airy sheet} is the natural two-parameter extension, and the \emph{directed landscape} is the result of the full, four-parameter limit. The directed landscape is a random function $\cL : \bR^4_\uparrow \to \bR$, where $\bR^4_\uparrow = \set{(x, t; y, s) \in \bR^4 : t < s \text{ or } (x, t) = (y, s)}$. The Airy sheet can be realised as the marginal $\cS(x, y) = \cL(x, 0; y, 1)$.

As we do not need the details of the constructions, we will content ourselves in recalling how they arise as limits of BLPP. Write $(x, t)_n = (t + 2 x n^{-1 / 3}, \floor{t n})$

\begin{theorem}[Theorem 1.5 of \cite{dov}]
\label{thm:BLPPDLLimit}
    Let $B = (B_1, B_2, \dots)$ be a sequence of standard Brownian motions, as before. The directed landscape appears as the limit
    \begin{equation*}
        n^{1 / 6}\brac[\big]{L_B \brac[\big]{(x, t)_n; (y, s)_n} - 2 (s - t)\sqrt{n} - 2 (y - x) n^{1 / 6}} \distto \cL(x, t; y, s).
    \end{equation*}
    The convergence here is in distribution, uniform on compact sets.
\end{theorem}

Both objects may be defined directly in terms of the Airy line ensemble. The details can be found in \cite{dauv-vir-lis}.

\subsection{Geodesics in the directed landscape}
Call the paths $\pi$ which attain the maximum in \eqref{eq:LPPDef} or \eqref{eq:BLPPTimes} or \eqref{eq:MelonDef} (with $k = 1$) \emph{geodesics}. When our variables are continuous these paths are unique almost surely, but to avoid ambiguity let us always take the right-most maximising path to be \emph{the} geodesic. 

We adopt a somewhat awkward convention for viewing these as non-decreasing cadlag functions, which will make our convergence result easier to state. Given points $a \le b \in \bZ^2$, the geodesic $\pi$ realising $L_\omega(a; b)$, and $v \in [a \cdot e_2, b \cdot e_2]$, define
\begin{equation*}
    \gamma_{\omega}(a; b; v) = \min\set{u \in \bZ: (u,\lfloor v\rfloor) \in \pi}.
\end{equation*}

Similarly, given $x \le y \in \bR$ and $m \le n \in \bZ$, let $\pi : [x, y] \to \set{m, \dots, n}$ be a path attaining the maximum in $\eqref{eq:BLPPTimes}$ (or with the collection of Brownian motions replaced by some line ensemble $F$), and take $v \in [m, n]$. Define
\begin{equation*}
    \gamma_{F}(x, m; y, n; v) = \min\set{u \in  [x, y]:  \pi(u) = \floor{v}}.
\end{equation*}
Essentially, the functions $\gamma$ encode the geodesics by taking a position on the vertical axis and giving the position on the horizontal axis at which the path reaches this level.

We can also define geodesics in the directed landscape. For $(x, t; y, s) \in \bR^4_\uparrow$, define $\Gamma(x, t; y, s; v)$ to be the almost surely unique path with the property that for each $v \in [t, s]$, one has
\begin{equation*}
    \cL(x, t; y, s) = \cL(x, t; \Gamma(x, t; y, s; v), v) + \cL(\Gamma(x, t; y, s; v), v; y, s).
\end{equation*}
There exists a great literature on the properties of these geodesics, but we content ourselves with an elementary distributional equivalence.

\begin{proposition}
\label{prop:DLGeoRescaling}
    Take $x < y$ and $t < s$. Then
    \begin{equation*}
        \set*{\frac{\Gamma\brac[\big]{x, t; y, s; t + v(s - t)} - x - v(y - x)}{(s - t)^{2 / 3}} : v \in [0, 1]} \disteq \set{\Gamma(0, 0; 0, 1; v) : v \in [0, 1]}.
    \end{equation*}
\end{proposition}
This is an application of Lemma 7 of \cite{bhat-duality} and the symmetries listed above it. We will simply write $\Gamma(v)$ for $\Gamma(0, 0; 0, 1; v)$.

\subsection{Modes of convergence}
We mention that various modes of convergence to the directed landscape are discussed in \cite{dauv-vir-lis}. Of these, the two relevant to us are \emph{compact}\footnote{In \cite{dauv-vir-lis}, they refer to \emph{graph} convergence, which is equivalent in this context.} and \emph{hypograph} convergence. By compact convergence we will mean distributional convergence with respect to uniform convergence on compact sets. Of hypograph convergence, it suffices for our purposes to know that it is weaker than compact convergence, and that it implies almost sure uniform convergence of geodesics. See Section 8 of \cite{dauv-vir-lis} for the precise statements.

\subsection{Strong approximations by Brownian motion}
Recall the environment of weights $\set{\omega(a) : a \in \bZ^2}$. Define a line ensemble $F = (F_1, F_2, \dots)$ with domain $\bR$ as follows:
\begin{itemize}
    \item When $u \ge 1$ is an integer, set $F_i(u) = \sum_{j = 0}^{u - 1}\omega(j, i)$.
    \item When $u \le 0$ is an integer, set $F_i(u) = \sum_{j = u}^{-1}\omega(j, i)$.
    \item For non-integer $u \in \bR$, let $F_i(u)$ be given by the linear interpolation between $F_i(\floor{u})$ and $F_i(\ceil{u})$.
\end{itemize}

When the $\omega$ have finite variance, Donsker's theorem tells us that the ensemble $F$ converges, after centring and rescaling, to a collection of independent Brownian motions. When additional moments are assumed, this convergence can be made quantitative by way of strong couplings to Brownian motions. 

The most famous result in this direction, and the one utilised in \cite{bodi-mart}, is due to  Koml\`os, Major and Tusn\`ady \cite{kmt-kmt,major-kmt}. While this would suffice for hypograph convergence, we are able to upgrade to compact convergence using the yet stronger coupling due to Sakhanenko \cite{Sakhanenko1985}. We use the translation found in \cite{Zaitsev2013}. Recall that we assume the $\omega$ have finite $p$-th moment, for $p > 2$.

\begin{proposition}[{\cite[Theorem 6]{Zaitsev2013}}]
    \label{prop:StrongApprox}
    Write $F = F_i$ for one of the independent random walks in our ensemble. There is a constant $c$ and a coupled Brownian motion $B$ such that for all $n$,
    \begin{equation}
        \label{eq:StrongApprox}
        \Ex{\max_{u \in [0, n]} \abs{F(u) - B(u)}^p} \le c n.
    \end{equation}
\end{proposition}

Strictly speaking, the original statement only covers integer times. We prove the (elementary) extension to all real times at the end of \cref{sec:proofs}. Jensen's inequality yields a useful corollary.

\begin{proposition}
    \label{prop:StrongApproxEx}
    In the same setting as \cref{prop:StrongApprox}, we have
    \begin{equation}
        \label{eq:StrongApproxEx}
        \Ex{\max_{u \in [0, n]} \abs{F(u) - B(u)}} \le c n^{1/p}.
    \end{equation}
\end{proposition}

\section{Results}
\label{sec:results}
We will establish versions of \cref{thm:BLPPALELimit,thm:BLPPDLLimit} for general last passage times, but under a different scaling. Rather than boxes of dimension $n \times n^\alpha$ with $0 < \alpha < 1$, as used in \cite{bodi-mart}, it will be notationally convenient to use dimensions $n^\beta \times n$, with $\beta > 1$. 

To allow the statements to fit horizontally on the page, let $(x, t)_{n, \beta} = (t n^{\beta} + 2 x n^{\beta - 1 / 3}, \floor{t n})$ and define
\begin{equation}
\label{eq:RescalledMelonDef}
    W^{n, \beta} F(y) = n^{(1 - 3 \beta)/ 6}W^n F\brac[\big]{0, 1; (y, 1)_{n, \beta}} - 2 n^{2 / 3} - 2 y n^{1 / 3}.
\end{equation}
Also define
\begin{equation}
\label{eq:RescalledLPPDef}
    d^{n, \beta}_\omega(x,t; y, s) = n^{(1 - 3 \beta)/ 6}L_\omega \brac[\big]{(x, t)_{n, \beta}; (y, s)_{n, \beta}} - 2 (s - t)n^{2 / 3} - 2 (y - x) n^{1 / 3},
\end{equation}
and similarly for $d^{n, \beta}_B$.

Our first result is Airy line ensemble convergence for general random walks under this scaling.
\begin{theorem}
\label{thm:GenALELimit}
We have $W^{n, \beta} F \distto \frakA$ compactly, whenever $\beta > \tfrac{7}{3}\cdot \tfrac{p}{p - 2}$.
\end{theorem}

From the Airy line ensemble limit we can infer a directed landscape limit for the passage times. This is more or less immediate given the tools of \cite{dauv-vir-lis}, except that there is an annoying technicality in the conversion between the fully discrete LPP and the semi-discrete BLPP.
\begin{theorem}
\label{thm:GenDLLimit}
We have $d^{n, \beta}_\omega \distto \cL$ whenever $\beta > \tfrac{7}{3}\cdot\tfrac{p}{p - 2}$, where this limit is in the sense of hypograph convergence. Moreover, this convergence is compact if $p > 5$. 
\end{theorem}

In particular, if $\omega$ has finite moments of all orders, then these hold with $\beta > 7 / 3$. The conjectured directed landscape limit for LPP corresponds to taking $\beta = 1$. 

That hypograph convergence implies uniform convergence of geodesics \cite[Theorem 8.7]{dauv-vir-lis} gives us the following distributional limits for geodesics near the axis.

\begin{corollary}
    Take $(x, t; y, s) \in \bR^4_\uparrow$. Let $\gamma_n = \gamma_\omega((x, t)_{n, \beta}; (y, s)_{n, \beta})$ be the cadlag function encoding a geodesic for $L_\omega$, taking the integer parts of the arguments if needed. For $\beta$ as above,
    % \begin{equation*}
    %     \set*{\frac{\gamma_n\brac[\big]{\floor{n t} + v(\floor{ns} - \floor{n t})} - n^{\beta}(t + v(s - t))}{2n^{\beta - 1 / 3}} : v \in [0, 1]} \distto \set{\Gamma(x, t; y, s; t + v(s - t)) : v \in [0, 1]}
    % \end{equation*}
    % and
    \begin{multline*}
        \set*{\frac{\gamma_n\brac[\big]{\floor{n t} + v(\floor{ns} - \floor{n t})} - n^{\beta}(t + v(s - t)) - 2 n^{\beta - 1 / 3}(x + v(y - x))}{2(s - t)^{2 / 3}n^{\beta - 1 / 3}} : v \in [0, 1]}\\ \distto \set{\Gamma(v) : v \in [0, 1]},
    \end{multline*}
    where the distributional convergence is with respect to uniform convergence.
\end{corollary}

\begin{remark}
    When we specialise to $(x, t; y, s) = (0, 0; 0, 1)$, this turns into the much more palatable
    \begin{equation*}
        \set*{\frac{\gamma_\omega\brac[\big]{0, 0; \floor{n^\beta}, n; v n} - v n^{\beta}}{2n^{\beta - 1 / 3}} : v \in [0, 1]} \distto \set{\Gamma(v) : v \in [0, 1]}.
    \end{equation*}
\end{remark}

These limits are just a combination of geodesic convergence and the normalisation found in \cref{prop:DLGeoRescaling}. From this we can identify the transversal fluctuation exponent in our stretched regime as $\beta - 1 / 3$. Up to a difference in parameterisation, this is the exponent predicted in \cite{bodi-mart}.

\section{Proofs}
\label{sec:proofs}
We will follow \cite{bodi-mart} and show that the strong coupling is strong enough to imply that the limits of \cref{thm:GenALELimit,thm:GenDLLimit} hold in general once they hold for Brownian motion. The limits in the case of Brownian motion are a simple application of Gaussian scaling: $L_B(x, t; y, s) \disteq a^{- 1 / 2}L_B(a x, t; a y, s)$. We arrive at the following.
\begin{lemma}
\label{lem:BLPPRescaledLimits}
    For all $0 \le \beta < \infty$, we have $W^{n, \beta}B \distto \frakA$ and $d^{n, \beta}_B \distto \cL$, both limits being in the sense of compact convergence.
\end{lemma}

With this, we proceed to prove the desired limit with general weights.

\begin{proof}[Proof of \cref{thm:GenALELimit}]
Couple $F$ and $B$ as in \cref{prop:StrongApprox}. Write
\begin{equation*}
    (T^k F)(y) = \sum_{i = 1}^k (W^{n}F)_i (y) = \sup_{\pi \in Q^k ((0, 1)_{n, \beta}; (y, n)_{n, \beta})} \int d F \circ \pi,
\end{equation*}
where for the second equality we recall the definition of the melon map. 

Recall that we view $W^{n,\beta}F$ and $\frakA$ as random functions from $ \mathbb{N} \times \mathbb{R}$ to $\mathbb{R}$. Thus, to show that $W^{n,\beta}F \xrightarrow{d} \frakA$ uniformly over compact sets, it suffices to consider the compact sets of the form $\{1, 2, \dots, \ell\} \times [w,z]$. Also, let us take $n$ large enough to ensure
\begin{equation*}
0 < n^\beta + 2 w n^{(3\beta - 1)/3} \le n^\beta + 2 z n^{(3\beta - 1)/3} < 2 n^\beta.
\end{equation*}
Then for $y \in [w, z]$ and $k \in \{1, \dots, \ell\}$, we have
\begin{align*}
    \abs{T^k F(y) - T^k B(y)} &\le T^k \abs{F - B}(y)\\
                              &= \sup_{\pi \in Q^k ((0, 1)_{n, \beta}; (y, 1)_{n, \beta})} \int d \abs{F - B} \circ \pi\\
                              &\le k \sup_{\pi \in Q^1 ((0, 1)_{n, \beta}; (y, 1)_{n, \beta})} \int d \abs{F - B} \circ \pi\\
                              &\le k \sum_{i = 1}^{n}\max_{u \in [0, n^\beta + 2y n^{\beta - 1 / 3}]}\abs{F_i(u) - B_i(u)}\\
                              &\le k \sum_{i = 1}^{n}\max_{u \in [0, 2 n^{\beta}]}\abs{F_i(u) - B_i(u)}.
\end{align*}
The last bound is uniform over all $y \in [w, z]$, so letting $\hat{T}^k[w, z] = \max_{y \in [w, z]}\abs{T^k F(y) - T^k B(y)}$, we have
\begin{equation*}
    \hat{T}^k[w, z] \le k \sum_{i = 1}^{n}\max_{u \in [0, 2 n^{\beta}]}\abs{F_i(u) - B_i(u)}.
\end{equation*}
Apply \cref{prop:StrongApproxEx} to this bound to find
\begin{equation*}
\label{eq:MelonBound}
    \Ex{\hat{T}^k[w, z]} \le c k n^{1 + \beta / p}.
\end{equation*}
Now Markov's inequality tells us that for any $\epsilon > 0$,
\begin{equation*}
    \Prob*{\max_{y \in [w, z]}\abs{T^k F(y) - T^k B(y)} > n^{1 + \beta / p + \epsilon}} \le c k n^{1 + \beta / p} n^{-({1 + \beta / p + \epsilon})} \to 0.
\end{equation*}

Given our assumption on $\beta$, we can choose $\epsilon$ such that $1+ \beta/p + \epsilon < (3\beta - 1)/6$. Then for all $k \in \{1,\dots, \ell\}$, we have
\begin{align*}
    \max_{y \in [w, z]} n^{(1-3\beta)/6}\abs{T^k F(y) - T^k B(y)} \rightarrow 0
\end{align*}
in probability. This implies that
\begin{align}\label{eq:converge_in_p}
    \max_{k \in \{1, \dots, \ell\}} \max_{y \in [w, z]} n^{(1-3\beta)/6}\abs{W^{n, \beta}_k F(y) - W^{n, \beta}_k B(y)} \rightarrow 0,
\end{align}
also in probability. 

To prove the claimed convergence to the Airy line ensemble, let $f : C(\set{1, \dots, \ell} \times [w, z]) \to [0, 1]$ be a uniformly continuous function, with its argument in the space of continuous functions on our compact set equipped with the $L^\infty$ norm. For any $\epsilon > 0$, we may find $\delta > 0$ such that if $G,\, H \in C(\set{1, \dots, \ell} \times [w, z])$ have $\norm{G - H}_\infty < \delta$, then $\abs{f(G) - f(H)} < \epsilon$. Thus, after restricting our melons to $\set{1, \dots, \ell} \times [w, z]$, we have the bound 
\begin{equation}
\label{eq:converge_in_d}
\abs[\big]{\Ex{f(W^{n,\beta}F)} -\Ex{f(W^{n,\beta}B)}} \leq \epsilon + \Prob{\norm{W^{n,\beta}F - W^{n,\beta}B}_\infty > \delta}
\end{equation} 
Due to the convergence in probability established in (\ref{eq:converge_in_p}), the right-hand side of (\ref{eq:converge_in_d}) can be made arbitrarily small. Combined with \cref{lem:BLPPRescaledLimits}, this gives the result.
\end{proof}

Having shown convergence of the melon, we would like to pass to a full directed landscape limit. In the semi-continuous setting of BLPP, the height of the top line of the melon is a last passage value, but this correspondence doesn't hold up exactly for the lattice on which the passage times $L_\omega$ live.

Consider $L_F$, defined analogously to $L_B$ in \eqref{eq:BLPPTimes}. Observe that $(W^n F(y))_1 = L_F(0, 1; y, n)$, as with BLPP. The top line of $W F$ also admits an interpretation as a lattice passage time.
\begin{lemma}
\label{lem:WFasPassageTime}
    For $a \le b \in \bZ^2$, we have the representation
    \begin{equation}
    \label{eq:WFasPassageTime}
        L_F(a; b) = \max_{\pi \in \Pi(a; b)}\sum_{\ov{r \in \pi}{r + e_1 \in \pi}} \omega(r).
    \end{equation}
    The maximum is over weakly increasing nearest-neighbour paths in $\bZ^2$ connecting $a$ to $b$, but now we only add weights from vertices from which we take a horizontal step.
    \begin{proof}
        Recall that we define
        \begin{equation}
        \label{eq:SWLPPTimes}
            L_F(x, m; y, n) = \max_{x = z_{m - 1} \le z_{m} \le \cdots \le z_{n} = y}\sum_{i = m}^{n} \sqrbrac{F_{i}(z_{i}) - F_{i}(z_{i - 1})}.
        \end{equation}
        There will always be an extremising path with its jumps occurring at extreme points of the line ensemble. As the lines of $F$ are linearly interpolated between integer times, these extreme points only exist at integer times. Thus when our end points are integers, we need only consider paths with integer jump points. Observe that a path in \eqref{eq:SWLPPTimes} collects weight only when it travels horizontally along a curve. That is, only when it moves horizontally on the lattice. This is equivalent to \eqref{eq:WFasPassageTime}.
    \end{proof}
\end{lemma}

This representation enables a simple comparison to $L_\omega$.
\begin{lemma}
    \label{lem:FMinusOmegaBound}
    When defined using the same weights, we have a deterministic inequality
    \begin{equation*}
        \abs[\big]{L_F(a; b) - L_\omega(a; b)} \le \max_{\pi \in \Pi(a; b)}\sum_{\ov{r \in \pi}{r + e_2 \in \pi}} \abs{\omega(r)}.
    \end{equation*} 
\end{lemma}

That is, the effect of disregarding weights on vertical steps is bounded by the maximum possible contribution of vertical steps, after taking absolute values. This is useful because the paths we consider take proportionally few vertical steps.

\begin{lemma}
    \label{lem:FMinusOmegaTail}
    Given $t > 0$, the following bound holds for the discrepancy between $L_F$ and $L_\omega$ on the box $[0, t n^\beta] \times [0, t n]$: 
    \begin{equation}
    \label{eq:FMinusOmegaTail}
        \Ex*{\max_{a,\, b \in ([0, t n^\beta] \times [0, t n]) \cap \bZ^2} \abs[\big]{L_F(a; b) - L_\omega(a; b)}} \le (tn)^{1 + \beta / p} \brac[\big]{\Exabs{\omega}^p}^{1 / p}.
    \end{equation}
    In particular, 
    \begin{equation}
    \label{eq:FMinusOmegaTailProb}
        \Prob*{\max_{a,\, b \in ([0, t n^\beta] \times [0, t n]) \cap \bZ^2} \abs[\big]{L_F(a; b) - L_\omega(a; b)} > n^{1 + \beta / p + \epsilon}} \to 0.
    \end{equation}
    \begin{proof}
        For simplicity take $t = 1$. Set $M_j = \max_{0 \le i \le n^\beta} \abs{\omega(i, j)}$. Then clearly from \cref{lem:FMinusOmegaBound}, the discrepancy is bounded with
        \begin{equation*}
            \abs[\big]{L_F(a; b) - L_\omega(a; b)} \le \sum_{j = a \cdot e_2}^{b \cdot e_2 - 1}M_j \le \sum_{j = 0}^{n - 1}M_j,
        \end{equation*}
        the latter bound uniform for all $a,\, b \in ([0, n^\beta] \times [0, n]) \cap \bZ^2$. Call the sum on the right $M$.

        The $p$-th moment of the $\omega$ is finite, thus we have the bound $\Exabs{M_j}^p \le n^\beta \Exabs{\omega}^p$. We then also get $\Exabs{M}^p < n^{p + \beta} \Exabs{\omega}^p$. Taking $p$-roots gives the desired bound. The claimed convergence of the probabilities is just Markov's inequality.
    \end{proof}
\end{lemma}

\begin{lemma}
\label{lem:FToOmegaCompact}
    When $\beta > \tfrac{7}{3}\cdot\tfrac{p}{p - 2}$, we have $\abs{d^{n, \beta}_\omega - d^{n, \beta}_F} \distto 0$ compactly.
    \begin{proof}
    Consider the difference
    \begin{equation*}
        \abs{d^{n, \beta}_F(x, t; y, s) - d^{n, \beta}_\omega(x, t; y, s)} = n^{(1 - 3 \beta) / 6}\abs{L_F((x, t)_{n, \beta}; (y, s)_{n, \beta}) - L_\omega((x, t)_{n, \beta}; (y, s)_{n, \beta})}.
    \end{equation*}
    By \cref{lem:FMinusOmegaTail}, the quantity in the absolute value is smaller than $n^{1 + \beta / p + \epsilon}$ with high probability, uniform for $(x, t; y, s)$ in some compact set. Our condition on $\beta$ is precisely enough for there to be $\epsilon > 0$ with $1 + \beta / p + \epsilon + (1 - 3 \beta) / 6 < 0$. The convergence follows.
    \end{proof}
\end{lemma}

Now we verify that the directed landscape limit holds for metrics $d^{n, \beta}_F$, defined analogously to \eqref{eq:RescalledLPPDef}. 
\begin{proposition}
\label{prop:SWLPPDLLimit}
    We have $d^{n, \beta}_F \distto \cL$ whenever $\beta > \tfrac{7}{3}\cdot\tfrac{p}{p - 2}$, where this limit is in the sense of hypograph convergence. This convergence is compact if $p > 5$.
\begin{proof}
    Theorems 1.14 and 1.16 of \cite{dauv-vir-lis} together give criteria for hypograph convergence, telling us that we have convergence to the directed landscape if there is suitable convergence of the melon to the Airy line ensemble, and if the one-point distributions match those of $\cL$. The former is exactly \cref{thm:GenALELimit}. For the latter, our argument for \cref{thm:GenALELimit} shows in particular that one-point distributions of the top lines of $W^{n, \beta} F$ and $W^{n, \beta} B$ agree in the limit. As recorded in \cref{lem:BLPPRescaledLimits}, the latter is known to converge to $\cL$.  

    There is a further criterion, Theorem 1.17 of \cite{dauv-vir-lis}, telling us when we can conclude compact convergence from hypograph convergence. It is sufficient that
    \begin{equation*}
    \label{eq:MomentCriterion}
        \limsup_{n \to \infty}\sup_{y \in [-3, 3]}\Exabs{d^{n, \beta}_F (0, 0; 0, y)^-}^{5 + \delta} < \infty 
    \end{equation*}
    for some $\delta > 0$. Here $d^-$ means the negative part of $d$. 
    
    Suppose $p > 5$ and let us choose $\delta > 0$ such that $q = 5 + \delta < p$ (although large enough to ensure $\beta > \tfrac{7}{3}\cdot\tfrac{q}{q - 2})$. Recall the proof of \cref{thm:GenALELimit} and specialise to $k = 1$. Recall also that the top lines of our melons record the passage times. We get for $n$ large enough that
    \begin{equation*}
        \max_{u \in [-3, 3]}\abs{L^{n, \beta}_F (0, 0; 0, u) - L^{n, \beta}_B (0, 0; 0, u)} \le \sum_{i = - 3 n}^{3 n}\max_{u \in [-6 n^{\beta}, 6 n^{\beta}]}\abs{F_i(u) - B_i(u)}.
    \end{equation*}
    Taking $q$-th moments and applying \cref{prop:StrongApprox} with Minkowski's inequality,
    \begin{equation*}
        \Ex{\max_{u \in [-3, 3]}\abs{L^{n, \beta}_F (0, 0; 0, u) - L^{n, \beta}_B (0, 0; 0, u)}^q} \le c (6 n)^q (6 n^\beta) = C n^{q + \beta}.
    \end{equation*}
    Now the rescaled passage times have
    \begin{equation*}
        \sup_{u \in [-3, 3]}\Exabs{d^{n, \beta}_F (0, 0; 0, u) - d^{n, \beta}_B (0, 0; 0, u)}^{q} \le C n^{q(1 - 3 \beta) / 6 + q + \beta}.
    \end{equation*}
    Our assumption on $\beta$ is such that this exponent is negative, and thus the quantity vanishes in the limit. We see that the condition of \eqref{eq:MomentCriterion} holds for $d^{n, \beta}_F$ if and only if it holds for $d^{n, \beta}_B$. But the condition is known to hold for $d^{n, \beta}_B$, owing to \cref{lem:BLPPRescaledLimits}.
\end{proof}
\end{proposition}

All that's left for \cref{thm:GenDLLimit} is to observe that as $\abs{d^{n, \beta}_\omega - d^{n, \beta}_F} \distto 0$ compactly (by \cref{lem:FToOmegaCompact}), and $d^{n, \beta}_F \distto \cL$ in either the hypograph sense or compactly (by \cref{prop:SWLPPDLLimit}), then --- the latter being a possibly weaker form of convergence --- we must have that $d^{n, \beta}_\omega \distto \cL$ in the same sense.

\begin{remark}
    We address a technical point which appears in the proof of Theorem 1.17 of \cite{dauv-vir-lis}, which we have used above. Strictly, the result is only proved for weights which are bounded below. To get around this, we define truncated weights $\tilde{\omega}_n = \omega \vee r_n$, where $r_n$ is a sequence diverging to $-\infty$ fast enough to ensure that $\Prob{\omega < r_n} = o(n^{-(1 + \beta)})$. This condition ensures that $\abs{d^{n, \beta}_{\omega} - d^{n, \beta}_{\tilde{\omega}_n}} \distto 0$ compactly. Thus $d^{n, \beta}_{\tilde{\omega}_n} \distto \cL$ in the hypograph sense.
    
    Also, it is easy to see that $(d^{n, \beta}_\omega)^-$ dominates $(d^{n, \beta}_{\tilde{\omega}_n})^-$. Thus the criterion in \eqref{eq:MomentCriterion} holds for $(d^{n, \beta}_{\tilde{\omega}_n})^-$ also. We then use Theorem 1.17 of \cite{dauv-vir-lis} to upgrade the convergence $d^{n, \beta}_{\tilde{\omega}_n} \distto \cL$. Using once more that $\abs{d^{n, \beta}_{\omega} - d^{n, \beta}_{\tilde{\omega}_n}} \distto 0$ compactly, now in the opposite direction, gives the desired limit.
\end{remark}

\begin{proof}[Proof of \cref{prop:StrongApprox}]
    Let $B$ be a coupled Brownian motion and write $\overline{B}$ for the process with $\overline{B}(n) = B(n)$ for integer $n$ and linear interpolation in between. Write
    \begin{equation*}
        \Delta[w, z] = \max_{u \in [w, z]}\abs{F(u) - B(u)},\quad \Delta_1[w, z] = \max_{u \in [w, z] \cap \bZ}\abs{F(u) - B(u)}.
    \end{equation*}
    We observe that $\Delta_1[w, z]$ may equally be written as $\max_{u \in [w, z]}\abs{F(u) - \overline{B}(u)}$, since this maximum must occur at an extreme point of the curves, namely an integer point. 
    
    In its original formulation, the approximation of Sakhanenko guarantees that we can find a coupled $B$ such that $\Ex{\Delta_1[0, n]^p} \le c n$. Now
    \begin{equation}
    \label{eq:StrongApproxTriangle}
        \Delta[w, z] \le \Delta_1[w, z] +  \max_{u \in [w, z]}\abs{B(u) - \overline{B}(u)}.
    \end{equation}
    Write $\Delta_2[w, z]$ for the latter term. Notice that $B(u) - \overline{B}(u)$ is just a concatenation of independent Brownian bridges. Write $R_m(t) = B(m + t) - \overline{B}(m + t)$ and $S_m = \max_{t \in [0, 1]} \abs{R_m(t)}$. In this notation, $\Delta_2[0, n] = \max_{0 \le m \le n - 1}S_m$. By the standard argument with the reflection principle, each $S_m$ is stochastically dominated by $\abs{Z^1_m} \vee \abs{Z^2_m}$, where $Z^1_m$, $Z^2_m$ are independent normals. So
    \begin{equation*}
        \Delta_2[0, n] \le \max_{0 \le m \le n - 1}(\abs{Z^1_m} \vee \abs{Z^2_m}).
    \end{equation*}
    Using the estimates for normal variables then, we can say $\Ex{\Delta_2[0, n]^p} \le C \log^{p / 2} n$. Returning to \eqref{eq:StrongApproxTriangle}, the statement follows from this and the bound of $\Delta_1[0, n]$. 
    
\end{proof}

\printbibliography
\end{document}